\documentclass[12pt]{amsart}
\usepackage[utf8]{inputenc}
\usepackage{amsmath,amssymb}

\begin{document}
\begin{abstract}
The generalization of Archimedes strategy to obtain the area
of a parabolic segment leads to combinatorial formulas involving minimal
covers of sets. These, in turn, are conjecturally related to $q$-binomial coefficients.
\end{abstract}

\title[Archimedes and minimal covers]{Archimedes' quadrature of the parabola and minimal covers}
\author{
Octavio Alberto Agustín Aquino
}
\address{Universidad de la Cañada, San Antonio Nanahuatipan Km 1.7 s/n.
Paraje Titlacuatitla, Teotitlán de Flores Magón, Oaxaca, México,
C.P. 68540.}
\email{octavioalberto@unca.edu.mx}
\address{Universidad Tecnológica de la Mixteca. Carretera a Acatlima Km 2.5 s/n. Huajuapan de León, Oaxaca, México, C. P. 69000.}
\email{octavioalberto@mixteco.utm.mx}
\date{February 16th, 2016}
\maketitle

\section{Introduction}

Archimedes \cite[pp. 233-252]{tH97} calculated the area any parabolic segment with supreme ingenuity via an exhaustion argument, successively
adding vertices to a polygon inscribed in it (beginning with a triangle).
Then he proved that each iteration contributed a fixed fraction
of the area of the polygon in the previous step (namely, one fourth of it), and thus he could sum the resulting geometric series.

When I planned to teach this rather elementary method in the unit interval, I realized that I could not use
the straightforward but convolved euclidean geometry, since my students were not familiar with it. Thus I recurred to the shoelace algorithm to compute the required
areas. This yields the determinant of a simple Vandermonde matrix to obtain the area of a single triangle in the unit interval with
$2^{n}$ equal subdivisions:
\[
\frac{1}{2}
\det\begin{pmatrix}
\frac{k-1}{2^{n}} & \left(\frac{k-1}{2^{n}}\right)^{2} & 1\\
\frac{k}{2^{n}} & \left(\frac{k}{2^{n}}\right)^{2} & 1 \\
\frac{k+1}{2^{n}} & \left(\frac{k+1}{2^{n}}\right)^{2} & 1
\end{pmatrix}
= \frac{1}{2}\left(\frac{2}{2^{3n}}\right) = \frac{1}{2^{3n}}.
\]

As it is readily seen, the value of the determinant is independent of the three consecutive points chosen to build the triangle.
Summing the area of the $2^{n}$ new triangles that appear in the $n$-th iteration, we obtain $\frac{2^{n}}{2^{3n}} = \frac{1}{2^{2n}} = \frac{1}{4^{n}}$,
as expected.

\section{Generalization}

What happens if we try the same trick with the function $f(x)=x^{s}$, with $s>2$? We are led to the determinant
\[
T(k,s,n):=
\det\begin{pmatrix}
\frac{k-1}{2^{n}} & \left(\frac{k-1}{2^{n}}\right)^{s} & 1\\
\frac{k}{2^{n}} & \left(\frac{k}{2^{n}}\right)^{s} & 1 \\
\frac{k+1}{2^{n}} & \left(\frac{k+1}{2^{n}}\right)^{s} & 1
\end{pmatrix} = \frac{(k+1)^{s}-2k^{s}+(k-1)^{s}}{2^{n(s+1)}}
\]

In this case, the area $\tfrac{1}{2}T(k,s,n)$ of the triangle depends on the points selected. For instance, when $s=3$, we have
\[
\frac{1}{2}\cdot\frac{(k+1)^{3}-2k^{3}+(k-1)^{3}}{2^{n(3+1)}} = \frac{3k}{2^{4n}}.
\]

To calculate, for example,
\[
\int_{0}^{1}x^{3}\,dx
\]
we first find the sum
\[
\sum_{i=0}^{2^{n-1}-1} \tfrac{1}{2}T(2i+1,3,n) = \sum_{i=0}^{2^{n-1}-1} \frac{3(2i+1)}{2^{4n}}
\]
that accounts for the area of the $n$-th iteration of the exhaustion. A virtue of this approach is that
only with the knowledge of the formula for triangular numbers $\sum_{k=0}^{n}k = \frac{n(n+1)}{2}$ we can find the area under
a cubic parabola in the unit interval, which is more difficult with the approach using Riemann sums with arithmetic or even
geometric subdivisions. More explicitly, we calculate
\[
\sum_{i=0}^{2^{n-1}-1} \frac{3(2i+1)}{2^{4n}} = \frac{1}{2^{4n+1}}\left(6\cdot\frac{2^{n-1}(2^{n-1}-1)}{2}+3\cdot 2^{n-1}\right)
= \frac{3}{2^{2n+2}}
\]
and thus
\[
\sum_{n=1}^{\infty} \frac{3}{2^{2n+2}} = \frac{3}{4}\left(\frac{1}{4(1-\frac{1}{4})}\right) = \frac{3}{4}\cdot \frac{1}{3} = \frac{1}{4}.
\]

Keeping in mind that we have calculated the area between the curve and the identity function, we simply subtract this result to $\frac{1}{2}$, and
we get the expected $\frac{1}{4}$. It is amusing to try by hand the case $s=4$.

It is worthwhile to mention that the former is a purely symbolic recast of Kirfel's more geometric exposition in \cite{cK13},
where upper bounds for the integrals are also considered. 

\section{A combinatorial excursion}

When $s\geq 2$, $T(2x+1,s,0)$ is a polynomial in $x$ of degree $s-2$.
If we arrange them in a triangle, we get
\[
\begin{matrix}
T(2x+1,2,0) = & 2,   &     &     &     &   \\
T(2x+1,3,0) = & 6+   & 12x,   &     &     &   \\
T(2x+1,4,0) = & 14+  & 48x+  & 48x^{2},   &     &   \\
T(2x+1,5,0) = & 30+  & 140x+ & 240x^{2}+  & 160x^{3},  &   \\
T(2x+1,6,0) = & 62+ & 360x+ & 840x^{2}+ & 960x^{3}+ & 480x^{4},
\end{matrix}
\]
whose coefficients can be rewritten as
\[
\begin{matrix}
1\cdot 2   &     &     &     &   \\
3\cdot  2   & 3\cdot 4   &     &     &   \\
7\cdot 2  & 12\cdot 4  & 6\cdot 8   &     &   \\
15\cdot 2  & 35\cdot 4 & 30\cdot 8  & 10\cdot 16  &   \\
31\cdot 2 & 90\cdot 4 & 105\cdot 8 & 60\cdot 16 & 15\cdot 32.
\end{matrix}
\]

Disregarding\footnote{The act
of overlooking the powers of two may seem artificial, and indeed it is done here
for the sake of the flow of ideas of the exposition. The combinatorial significance of the 
result was found as an accident, caused precisely by failing to consider only the triangles with odd indices during
the limit process.} the powers of $2$ that we factorized, we get $\tfrac{1}{2}T(x+1,s,0)$. Now
note the following
\begin{align*}
T(x+1,s,0) &= (x+2)^{s}-2(x+1)^{s}+x^{s}\\
 &= \sum_{k=0}^{s}\binom{s}{k}2^{s-k}x^{i}-2\sum_{k=0}^{s}\binom{s}{k}x^{k}+x^{s}\\
&=2\sum_{k=0}^{s}\binom{s}{k}(2^{s-k-1}-1)x^{k}+x^{s}\\
&=2\sum_{k=0}^{s}\binom{s}{k}(2^{k-1}-1)x^{s-k}+x^{s}\\
&= 2\sum_{k=2}^{s}\binom{s}{k}(2^{k-1}-1)x^{s}
\end{align*}

The coefficients $\binom{s}{k}(2^{k-1}-1)$ of the polynomial correspond to the
number $M(s,2,k)$ of minimal\footnote{
A \emph{minimal} cover is a cover such that the elimination of any of
its members results in a family of sets that fails to cover the original set.} $2$-covers of a labeled
$s$-set that cover $k$ points uniquely (where $k\geq 2$), listed as sequence A057963 in the OEIS \cite{oeis}. Hence
\[
T(x+1,s,0) = 2\sum_{k=2}^{s} M(s,2,k)x^{s-k}
\]
and
\[
T(x+1,s,1) = \frac{1}{2^{s}}\sum_{k=2}^{s} M(s,2,k)x^{s-k}.
\]

In general
\[
T(x+1,s,n) = \frac{1}{2^{n(s+1)}}\sum_{k=2}^{s}M(s,2,k)x^{s-k},
\]
and therefore
\begin{align}
\frac{s-1}{s+1}&=1-\frac{2}{s+1}\notag\\
&=1-2\int_{0}^{1}x^{s}\,dx= \sum_{n=1}^{\infty}\sum_{x=1}^{2^{n-1}-1}T(2x+1,s,n)\notag\\
&= \sum_{n=1}^{\infty}\frac{1}{2^{(n-1)(s+1)}}\sum_{x=1}^{2^{n-1}-1}\sum_{k=2}^{s} \frac{M(s,2,k)}{2^{k}}x^{s-k}.\label{E:SumaCubiertasMinimas}
\end{align}

Hearne and Wagner proved \cite{HW73} that if $M(s,j,k)$ is the number of
$j$-member minimal covers of an $s$-set that cover $k$ elements
uniquely, then
\[
M(s,j,k) = \binom{s}{k}(2^{j}-j-1)^{s-k}S(k,j)
\]
where $S(k,j)$ denotes the Stirling numbers of the second kind. This formula allows us to extend \eqref{E:SumaCubiertasMinimas}
for $j> 2$, and the following Maxima code (for the particular case when $s=4$ and $j=2$, which should
be changed accordingly), valid for $3\leq j \leq s-1$, automatizes the calculations.

\begin{verbatim}
sumando: subst([s=4,j=2],(1/2^((n-1)*(s+1)))*
 sum(sum(binomial(s,k)*(2^j-j-1)^(s-k)*stirling2(k,j)*
 x^(s-k)/2^k,k,2,s),x,0,2^(n-1)-1)), simpsum;
sumando: ratexpand(sumando);
acumulado:0;
for k:1 thru nterms(sumando)step 1 do
 acumulado:acumulado+sum(part(sumando,k),n,1,inf),simpsum;
print(acumulado);
\end{verbatim}

We can now compile the following triangle of rational numbers.
\begin{equation}\label{E:TrianguloGen}
\begin{matrix}
s\backslash j & 1 & 2 & 3 & 4 & 5 & 6 & 7 & 8\\[0.5em]
2 & \frac{1}{\bf{3}}  & \frac{1}{\bf{3}}  & & & & & & \\[0.5em]
3 & \frac{1}{\bf{7}}  & \frac{1}{2} & \frac{1}{7} & & & & &\\[0.5em]
4 & \frac{1}{\bf{15}} & \frac{3}{5} & \frac{10}{21}            & \frac{1}{15} & & & &\\[0.5em]
5 & \frac{1}{\bf{31}} & \frac{2}{3} & \frac{865}{651}          & \frac{71}{186}       &\frac{1}{31}            & & & &\\[0.5em]
6 & \frac{1}{\bf{63}} & \frac{5}{7} & \frac{2630}{651}         & \frac{1427}{\bf{651}}&\frac{181}{651}         &\frac{1}{63}              & & &\\[0.5em]
7 & \frac{1}{\bf{127}}& \frac{3}{4} & \frac{163133}{\bf{11811}}& \frac{306553}{15748} &\frac{36667}{\bf{11811}}&\frac{145}{762}           &\frac{1}{127}&\\[0.5em]
8 & \frac{1}{\bf{255}}& \frac{7}{9} & \frac{3368938}{66929}    &\frac{129115655}{602361}&\frac{46958822}{602361}&\frac{43662}{\bf{10795}} &\frac{4036}{32385} & \frac{1}{255}
\end{matrix}
\end{equation}

Using the notation of Hearne and Wagner $M(s,k) = \sum_{j=0}^{k}M(s,j,k)$
and their formula for the generating function
\begin{align*}
 M_{s}(x) &= \sum_{k=0}^{s}M(s,k)x^{k} = \sum_{k=0}^{s}\sum_{j=0}^{k}M(s,j,k)x^{k}\\
&= \sum_{j=0}^{s}\frac{1}{j!}\sum_{\ell=0}^{j}(-1)^{j-\ell}\binom{j}{\ell}(2^{j}-j-1+\ell x)^{s},
\end{align*}
where $M(s,j,k)$ are defined as zero where convenient, we have
\[
x^{s}M_{s}(\tfrac{1}{2x}) = \sum_{k=0}^{s}\sum_{j=0}^{s}\frac{M(s,j,k)}{2^{k}}x^{s-k}
\]
and thus we get 
\begin{multline*}
\sum_{n=1}^{\infty}\frac{1}{2^{(n-1)(s+1)}}\sum_{x=1}^{2^{n-1}-1}\sum_{k=2}^{s}\sum_{j=0}^{s}\frac{M(s,j,k)}{2^{k}}x^{s-k}\\
=\sum_{n=1}^{\infty}\frac{1}{2^{(n-1)(s+1)}}\sum_{x=1}^{2^{n-1}-1}x^{s}M_{s}(\tfrac{1}{2x})\\
=\sum_{n=1}^{\infty}\frac{1}{2^{(n-1)(s+1)}}\sum_{x=1}^{2^{n-1}-1} \sum_{j=0}^{s}\frac{1}{j!}\sum_{\ell=0}^{j}(-1)^{j-\ell}\binom{j}{\ell}((2^{j}-j-1)x+\tfrac{\ell}{2} )^{s}\\
=\sum_{n=1}^{\infty}\frac{1}{2^{n-1}}\sum_{x=1}^{2^{n-1}-1} \sum_{j=0}^{s}\frac{1}{j!}\sum_{\ell=0}^{j}(-1)^{j-\ell}\binom{j}{\ell}((2^{j}-j-1)\tfrac{x}{2^{n-1}}+\tfrac{\ell}{2} )^{s}
\end{multline*}
whence we can deduce the following lower bound for the sums of the rows of the triangle
\begin{multline*}
\sum_{n=1}^{\infty}\frac{1}{2^{(n-1)(s+1)}}\sum_{x=1}^{2^{n-1}-1}\sum_{k=2}^{s}\sum_{j=0}^{s}\frac{M(s,j,k)}{2^{k}}x^{s-k}\\
\geq \int_{0}^{1} \sum_{j=0}^{s}\frac{1}{j!}\sum_{\ell=0}^{j}(-1)^{j-\ell}\binom{j}{\ell}((2^{j}-j-1)x+\tfrac{\ell}{2} )^{s}\, dx\\
= \sum_{j=0}^{s}\frac{1}{j!}\sum_{\ell=0}^{j}(-1)^{j-\ell}\binom{j}{\ell}\left.\frac{((2^{j}-j-1)x+\tfrac{\ell}{2})^{s+1}}{(s+1)(2^{j}-j-1)}\right|_{0}^{1}\\
= \sum_{j=0}^{s}\frac{1}{j!}\sum_{\ell=0}^{j}(-1)^{j-\ell}\binom{j}{\ell}\frac{1}{s+1}\sum_{t=1}^{s+1}(2^{j}-j-1+\tfrac{\ell}{2})^{t-1}(\tfrac{\ell}{2})^{s-t+1}
\end{multline*}

\section{A conjecture}

The boldface denominators in \eqref{E:TrianguloGen} are the values when $q=2$ of the $q$-analogues of the binomial coefficients
\[
 \binom{s}{m}_{q} = \frac{[s]_{q}!}{[m]_{q}![s-m]_{q}!},
\]
where $[n]_{q} = \frac{q^{n}-1}{q-1}$ and the $q$-analogue of the factorial is defined inductively by $[0]_{q}!=1$ and $[n]_{q}!=[n]_{q}([n-1]_{q}!)$
 (see \cite[sequence A022166]{oeis}).
It is to be noted that some entries are divisors of the corresponding $2$-binomial coefficient. This pattern persist for
further rows of the triangle, but so far no explanation is evident to the author. The numerators are even more mysterious.
The sequences of numerators and denominators of the triangle are now entries A280752 and A280753 (respectively) in the OEIS.

\section*{Acknowledgement}

This note could not have been written without the valuable help of the OEIS. The author also thanks Michel Marcus for
pointing out an error regarding the ID number of one sequence in the OEIS in a previous version of this paper, and for his notification of the inclusion of the triangle obtained here in the encyclopedia.

\bibliographystyle{abbrv}
\bibliography{archimedes_integral}

\end{document}